\documentclass{article}
\usepackage{latexsym}
\usepackage{amsmath}
\usepackage{amsfonts}
\usepackage{euler}
\usepackage{graphicx}
\usepackage{hyperref}

\DeclareMathSymbol\leqslant{\mathrel}{AMSa}{"36}
\DeclareMathSymbol\geqslant{\mathrel}{AMSa}{"3E}
\DeclareMathSymbol\nmid{\mathrel}{AMSb}{"2D}
\renewcommand\le\leqslant
\renewcommand\ge\geqslant

\newtheorem{prop}{Proposition}[section]
\newtheorem{conj}[prop]{Conjecture}
\newtheorem{assertion}[prop]{Assertion}

\newenvironment{proof}{\par\noindent{\bf Proof: }}{\par}
\newenvironment{rem}{\par\noindent{\bf Remark: }}{\par}

\newcommand{\qed}{\hspace*{\fill}$\Box$}

\title{On Heron Simplices and Integer Embedding}
\author{Jan Fricke}

\begin{document}
\DeclareGraphicsExtensions{.eps}

\maketitle

\begin{abstract}
  In \cite{UPNT} problem D22 Richard Guy asked for the existence of simplices
  with integer lengths, areas, volumes\dots. In dimension two this is well
  known, these triangles are called {\em Heron triangles}. Here I will present
  my results on Heron tetrahedra, their connection to the existence of an
  integer box (problem D18), the tools for the search for higher dimensional
  Heron simplices and my nice embedding conjecture about Heron simplices,
  which I can only proof in dimension two, but I verified it for a large range
  in dimension three.
\end{abstract}

\section{Heron triangles}
\subsection{Basics}
\begin{prop} (Heron formula)
Let $a$, $b$ and $c$ be the lengths of the sides of a triangle, $s$ the
half of the perimeter and $A$ the area. Then holds
\begin{equation}\label{eq:Heron}
A^2=s(s-a)(s-b)(s-c).
\end{equation}
\qed
\end{prop}

Now we define: a {\em Heron triangle} is a triangle with integer sides and area.

\begin{prop}\label{prop:peri}
The perimeter of a Heron triangle is even.
\end{prop}
\begin{proof}
Assume the perimeter is odd. Then all factors on the right hand side of
Equation (\ref{eq:Heron}) are halves of odd integers, so the product is an odd
integer divided by $16$, which cannot be an integer. \qed
\end{proof}

\begin{prop}
  If the sides of a Heron triangle have a common factor $t$, then the area is
  divisible by $t^2$, i.e. the triangle scaled down by factor $t$ is also
  Heron. Furthermore the area of a Heron triangle is always divisible by $6$.
\end{prop}
\begin{proof}
It is sufficient to show this for $t$ prime. For $t\ne2$ it is obvious, for
$t=2$ consider the term $s(s-a)(s-b)(s-c)$ modulo 4 and one gets that it can't
be congruent $1$. The second part can in a similar way be checked modulo $4$
and $3$. \qed
\end{proof}

\subsection{Special types}
Now I will give some examples for Heron triangles. There are two main
classifications for triangles: \{acute, right-angled and obtuse\} and
\{regular, isosceles and generic\}. So we first search for right-angled and
isosceles Heron triangles.

\subsubsection{Right-angled triangles}
Right-angled triangles have to be Pythagorean but the reverse is also true:
\begin{prop}
Any Pythagorean triangle is a Heron triangle.
\end{prop}
\begin{proof}
The lengths of the sides are integer, so we only have to proof that the area
is integer. For the area $A$ holds $A=ab/2$, where $a$ and $b$ are the lengths
of the catheti. Hence if $A$ is not an integer then $a$ and $b$ have to be
odd. But then follows $a^2\equiv b^2 \equiv 1 \bmod 4$ and $c^2 = a^2+b^2\equiv
2 \bmod 4$ and that is not possible.\qed
\end{proof}

\subsubsection{Isosceles triangles}
\begin{prop}
Any Heron isosceles triangle is divided by the center line into two congruent
Pythagorean triangles.
\end{prop}
\begin{proof}
Using Proposition \ref{prop:peri} we get that the basis must be an even number,
say $2a$. Let $b$ the length of the arms and $h$ be the length of
the center line. Now we only have to proof that $h$ is an integer. We get
$s=(2a+b+b)/2=a+b$ and from the Heron formula
$A^2=s(s-2a)(s-b)(s-b)=(b^2-a^2)a^2$. From $A$ is an integer follows that
$b^2-a^2$ is a square number, but that is $h^2$, so $h$ is integer.\qed
\end{proof}

\begin{rem}
  There is not intersection of this two classes -- right-angled isosceles
  triangles don't have integer sides.
\end{rem}

\subsection{Embedding properties}
The examples in the previous section show that some classes of
Heron triangles can trivially be represented as lattice triangles.
First we show the following weaker result.
\begin{prop}
Any Heron triangle can be represented with rational coordinates.
\end{prop}
\begin{proof}
We only have to show that the height is rational and the length of
the height-tiles. But this follows from the equations
$h_c=\frac{A}{2c}$ and $c_a=\frac{-a^2+b^2+c^2}{2c}$.\qed
\end{proof}

There are Heron triangles without integer heights as the example
$(5,29,30)$ with area $72$ shows.

\begin{prop}
Any Heron triangle can be represented with integer coordinates.
\end{prop}
\begin{proof}
We know that it is representable with rational coordinates. Then
we can apply Proposition \ref{prop:em2}.\qed
\end{proof}

\section{Heron tetrahedra}
\subsection{Basics}
\begin{prop}
The volume $V$ of the tetrahedron $ABCD$ is given by
\begin{equation}\label{eq:v3}
V^2=\frac{1}{288}\left|
\begin{array}{ccccc}
 0   & AB^2 & AC^2 & AD^2 & 1 \\
BA^2 &  0   & BC^2 & BD^2 & 1 \\
CA^2 & CB^2 &  0   & CD^2 & 1 \\
DA^2 & DB^2 & DC^2 &  0   & 1 \\
 1   &  1   &  1   &  1   & 0
\end{array}
\right|.
\end{equation}
\end{prop}
\begin{proof}
This formula can be found in a not so intuitive notation e.g. in \cite{TB}.\qed
\end{proof}
The Equation (\ref{eq:v3}) gives also a criterion whether an
tetrahedron exists with the given edges:
\begin{prop}
A tetrahedron with given edges exists if and only if all faces
exist (i.e. in all triangles the triangle inequality holds, which
is equivalent to: the right hand side of Equation (\ref{eq:Heron}) is
positive) and the right hand side of Equation (\ref{eq:v3}) is
positive.\qed
\end{prop}

Now we define: a {\em Heron tetrahedron} is a tetrahedron where the lengths of
the edges, the area of the faces and the volume are integers.

\begin{prop}
  If the sides of a Heron tetrahedron have a common factor $t$, then the volume is
  divisible by $t^3$, i.e. the tetrahedron scaled down by factor $t$ is also
  Heron. Furthermore the volume of a Heron tetrahedron is always divisible by
  $336=2^4\cdot3\cdot7$.
\end{prop}
\begin{proof}
As for triangles it is sufficient to show this for $t$ prime. For $t\ne2,3$ it
is obvious. For $t=2,3$ one have to proof that if the determinant on the right
hand side of Equation (\ref{eq:v3}) is the double of a square and all areas of
the faces are integer (i.e. the right hand side of Equation (\ref{eq:Heron})
is a square) then it is divisible by $288$. Because I have no elegant proof
for that this was done by a ``brute force'' calculation modulo $9$ and modulo
$16$ using a computer program.

The second part was done in a similar way by checking modulo $2^{12}$, $3^3$
and $7$ using a computer. \qed
\end{proof}

\subsection{Special types}
Like in the case of triangles there are some special types of tetrahedra.

\subsubsection{Semi-regular or isosceles tetrahedra}
\begin{figure}[htb]
  \newsavebox{\srinabox}
  \newlength{\srinaboxlength}
  \sbox{\srinabox}{\includegraphics[scale=0.5]{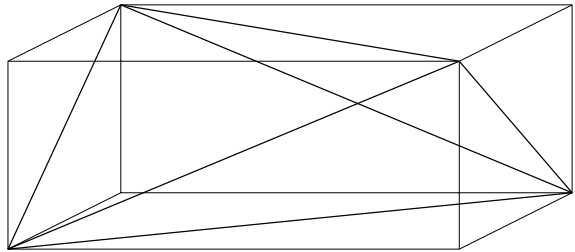}}
  \settowidth{\srinaboxlength}{\usebox{\srinabox}}
  \begin{minipage}[c]{\srinaboxlength}
    \includegraphics[scale=0.5]{srinbox}
  \end{minipage}
  \hfill
  \begin{minipage}[c]{3cm}
    \setlength{\unitlength}{1mm}
    \begin{picture}(30,40)
      \put(25,5){\makebox(0,0)[cc]{\ref{sr:circ}}}
      \put(25,15){\makebox(0,0)[cc]{\ref{sr:angle}}}
      \put(25,25){\makebox(0,0)[cc]{\ref{sr:net}}}
      \put(25,35){\makebox(0,0)[cc]{\ref{sr:peri}}}
      \put(15,5){\makebox(0,0)[cc]{\ref{sr:3p}{\tiny ci}}}
      \put(15,15){\makebox(0,0)[cc]{\ref{sr:3p}{\tiny mc}}}
      \put(15,25){\makebox(0,0)[cc]{\ref{sr:opp}}}
      \put(15,35){\makebox(0,0)[cc]{\ref{sr:cong}}}
      \put(5,5){\makebox(0,0)[cc]{\ref{sr:box}}}
      \put(5,15){\makebox(0,0)[cc]{\ref{sr:3p}{\tiny mi}}}
      \put(5,25){\makebox(0,0)[cc]{\ref{sr:hei}}}
      \put(5,35){\makebox(0,0)[cc]{\ref{sr:area}}}
      \put(15,27){\vector(0,1){6}}
      \put(5,17){\vector(0,1){6}}
      \put(15,17){\vector(0,1){6}}
      \put(25,17){\vector(0,1){6}}
      \put(5,7){\vector(0,1){6}}
      \put(25,7){\vector(0,1){6}}
      \put(5,30){\vector(0,1){3}}
      \put(5,30){\vector(0,-1){3}}
      \put(16,35){\vector(1,0){8}}
      \put(6,5){\vector(1,0){7}}
      \put(17,5){\vector(1,0){7}}
      \put(14,35){\vector(-1,0){8}}
      \put(7,7){\vector(1,1){6}}
      \put(7,33){\vector(1,-1){6}}
      \put(23,27){\vector(-1,1){6}}
      \put(23,33){\vector(-1,-1){6}}
      \put(10,15){\vector(1, 2){4}}
      \put(10,15){\vector(-1, -2){4}}
    \end{picture}
  \end{minipage}\\
  \begin{minipage}[t]{\srinaboxlength}
    \caption{Semi-regular tetrahedron in a box.}
    \label{fig:srinabox}
  \end{minipage}
  \hfill
  \begin{minipage}[t]{3cm}
    \caption{Proofs.}
    \label{fig:srproof}
  \end{minipage}
\end{figure}

For a tetrahedron the following properties are equivalent:
\begin{enumerate}
\itemsep0cm
\item\label{sr:opp} The opposite edges have the same length.
\item\label{sr:cong} All faces are congruent (acute) triangles.
\item\label{sr:peri} The perimeters of the faces coincide.
\item\label{sr:area} All faces have the same area.
\item\label{sr:hei} All body heights have the same length.
\item\label{sr:box} The tetrahedron can be inscribed in a box
  (see Figure \ref{fig:srinabox}).
\item\label{sr:3p} Two of the following three points coincide: the mass
  center, the center of the circumscribed and inscribed sphere.
\item\label{sr:circ} The circumcircle radii of the faces coincide.
\item\label{sr:angle} In all vertices the sum of the face angles equals $\pi$.
\item\label{sr:net} The net of the tetrahedron is a triangle with connected
  center points.
\end{enumerate}
\dots and much more. For a scheme of the proofs see Figure \ref{fig:srproof}.
If one of these properties is fulfilled we call the tetrahedron {\em semi-regular}. 
I also found the name {\em isosceles} for that (e.g. on \cite{EW}, there are a
lot of references), but I prefer the name semi-regular because these
tetrahedra are ``more regular than isosceles'', at least I will show in
Section \ref{sec:if} a generalization of isosceles triangles which fits better in this
context.

Hence to get a Heron semi-regular tetrahedron one has only to
take an acute Heron triangle and verify whether the volume is
integer, and the volume can be calculated by
$V^2=(a^2+b^2-c^2)(a^2-b^2+c^2)(-a^2+b^2+c^2)/72$.
So the smallest examples of this type can be easy computed as $(203,195,148)$,
$(888,875,533)$ and $(1804,1479,1183)$ its multiples. These were the first
known Heron tetrahedra. See Appendix \ref{ap:sr} for more examples of Heron 
semi-regular tetrahedra.

\subsubsection{Right-angled-vertex tetrahedra}
\begin{figure}[htb]
  \begin{center}
    \includegraphics[scale=0.5]{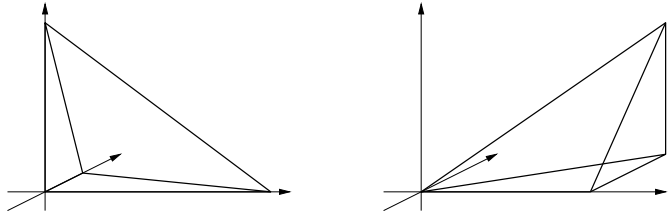}
    \caption{Right-angled-vertex and right-angled-face tetrahedron}
    \label{fig:ra}
  \end{center}
\end{figure}

The right-angled triangle can be generalized in two ways. The first type I will
call {\em right-angled-vertex tetrahedron}: a
tetrahedron where in one vertex $O$ all angles are right angles (see left
picture of Figure \ref{fig:ra}). Let $x$, $y$ and
$z$ be the lengths of these edges which ends in $O$. Then the area of the
faces which contains $O$ are automatically integer, and also the volume is
integer. So there's only left to examine the lengths $\sqrt{x^2+y^2}$,
$\sqrt{x^2+z^2}$ and $\sqrt{y^2+z^2}$, and the area of the face opposite to $O$.
This area can easy be computed by Formula (\ref{eq:Heron}) as
$\sqrt{x^2y^2+x^2z^2+y^2z^2}/2$.

\begin{prop}
There is a Heron right-angled-vertex tetrahedron if and only if there is an integer
box\footnote{An integer box is a box with integer edges, face diagonals and
body diagonal. See \cite{UPNT} problem D18.}.
\end{prop}

\begin{proof}
Let $u=yz$, $v=xz$ and $w=xy$ be the lengths of the box. Then the
lengths of the face diagonals are $x\sqrt{y^2+z^2}$, $y\sqrt{x^2+z^2}$
and $z\sqrt{x^2+y^2}$, and the length of the body diagonal is
$\sqrt{x^2y^2+x^2z^2+y^2z^2}$. So we have an integer box.

Lets $u$, $v$ and $w$ the lengths of an integer box. Set $x=vw$,
$y=uw$ and $z=uv$. Then the lengths of the other three edges
are $u\sqrt{v^2+w^2}$, $v\sqrt{u^2+w^2}$ and $w\sqrt{u^2+v^2}$,
which are integers. The area of the to $O$ opposite face is
$uvw\sqrt{u^2+v^2+w^2}/2$ which is also an integer.\qed
\end{proof}

But: the problem of the existence of an integer box is an ``notorious unsolved
problem''(\cite{UPNT}, problem D18).

\subsubsection{Right-angled-face tetrahedra}
The second type I will call {\em right-angled-face tetrahedron}: a tetrahedron
whose vertices $A_0$, $A_1$, $A_2$, $A_3$ can be represented by the Cartesian
coordinates $A_0=(0,0,0)$, $A_1=(a_1,0,0)$, $A_2=(a_1,a_2,0)$,
$A_3=(a_1,a_2,a_3)$ (see right picture of Figure \ref{fig:ra}).
Note that this condition is equivalent to all faces 
being right-angled triangles. Let's examine under which conditions this will
generate a Heron tetrahedron. First all edges have to be integer,
i.e. $a_1$, $a_2$, $a_3$, $\sqrt{a_1^2+a_2^2}$, $\sqrt{a_1^2+a_2^2+a_3^3}$
and $\sqrt{a_2^2+a_3^2}$. The faces are Pythagorean triangles, so their area is
integer and the volume is $a_1a_2a_3/6$ which is an integer (the product of
the catheti is always divisible by $6$). Now take a look to \cite{UPNT} problem
D18 and you get to know that this problem is equivalent to find a nearly
integer box, where only one edge is not integer. The smallest
solutions $(672,104,153)$ and $(756,117,520)$ were known to Euler. Hence these
solutions gives us the two smallest right-angled-face tetrahedra. See Appendix
\ref{ap:raf} for more.

\subsubsection{Isosceles-face tetrahedra}\label{sec:if}
Like in the case of right-angled-face tetrahedra one can ask whether there are
Heron tetrahedra which have only isosceles triangles as faces. One can find a
combinatorical classification:
\begin{enumerate}
\item $\overline{AB}=\overline{AC}=\overline{AD}$
(such Heron tetrahedra seem not to exist),
\item $\overline{AB}=\overline{BC}=\overline{CD}$ 
and $\overline{CA}=\overline{AD}=\overline{DB}$ (such Heron tetrahedra don't
exist: assume there is one, then from the classification of Heron isosceles
triangles follows, that the two lengths are doubles of integer $a$ and $b$
with $a^2+u^2=(2b)^2$ and $b^2+v^2=(2a)^2$. Adding this equations one gets
$u^2+v^2=3(a^2+b^2)$ which implies $a=b=u=v=0$) and
\item $\overline{AB}=\overline{BC}=\overline{CD}=\overline{DA}$.
\end{enumerate}
The last type I call {\em isosceles-face} tetrahedra. The nice analog to the
two-dimensional case is, that:
\begin{prop}
  Any Heron isosceles-face tetrahedron is divided by the center line into four
  congruent Heron right-angled-face tetrahedra.
\end{prop}
Here the center line is the line which connects the center of $AC$ with the
center of $BD$.

\begin{figure}[htb]
  \begin{center}
    \includegraphics[scale=0.5]{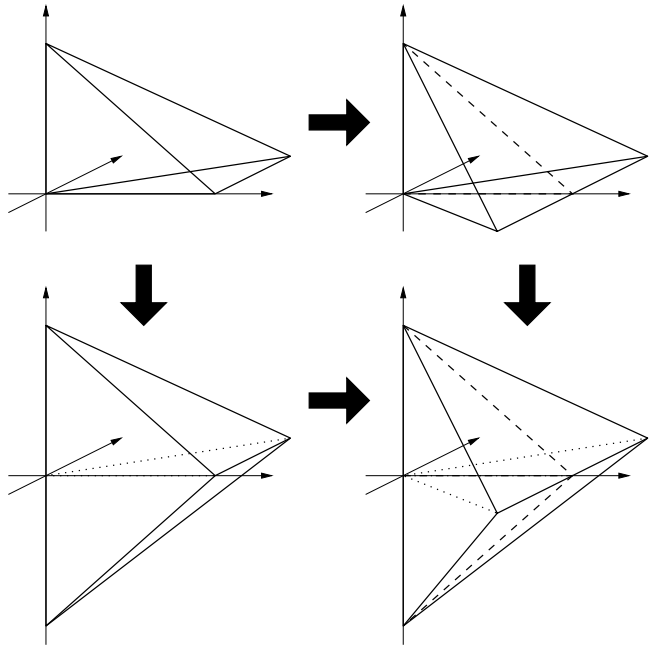}
    \caption{Construction of an isosceles-face tetrahedron from a
      right-angled-face tetrahedron}
    \label{fig:ra2if}
  \end{center}
\end{figure}

Figure \ref{fig:ra2if} shows how to construct an isosceles-face tetrahedron from a 
right-angled-face tetrahedron and the decomposition of an isosceles-face
tetrahedron into four right-angled-face tetrahedron (picture on the lower right).

\begin{proof}
  $AC$ and $BD$ are bases of isosceles triangles, so they have lengths $2a$
  and $2b$ respectively. Let $x$ be the length of the other four edges. We
  have to show that the length $m$ of the center line is integer. From
  $V=4\cdot(\frac16abm)$ follows that $m$ is rational and from
  $a^2+m^2+b^2=x^2$ follows that $m^2$ is integer. Hence $m$ is integer.\qed
\end{proof}

See Appendix \ref{ap:raf} for examples.

\subsubsection{Intersections and compositions}\label{sec:othertypes}
Now we look for the intersection of this classes:

\noindent
\begin{tabular}{|l|c|c|c|c|}
  \hline
  & r.-a.-v. & r.-a.-f. & s.-r. & i.-f. \\
  \hline
  right-angled-vertex & $\times$ & g & g & r-i \\
  \hline
  right-angled-face & g & $\times$ & g & r-i \\
  \hline
  semi-regular & g & g & $\times$ & o \\
  \hline
  isosceles-face & r-i & r-i & o & $\times$ \\
  \hline
\end{tabular}

\noindent Here ``g'' means ``not possible from geometric reasons'' (remember
that the faces of a semi-regular tetrahedron are acute), ``r-i'' means ``not
possible because there is no Heron right-angled isosceles triangle'' and ``o''
stands for an open case. The existence of such tetrahedra is equivalent to the
existence of non-trivial integer solutions of $x^4+y^4=2z^2$ (i.e. other than
$x^2=y^2=|z|$).

So these classes of Heron tetrahedra seem to be disjoint.

Finally let me make some remarks on other types of tetrahedra which may play a
role in this context. If one looks at the table in Appendix \ref{ap:gen} one
will notice that there some tetrahedra which have two isosceles triangles with
common basis as faces. They are composed of two congruent tetrahedra where one
edge is perpendicular to one face. (The construction is quite the same as
shown in Figure \ref{fig:ra2if}.) Using the embedding property of Heron
triangles it's easy to see, that these types can be represented with integer
coordinates, but this will be the topic of the next section.

\subsection{Embedding properties}
As in the case of Heron triangles the following weak result holds.
\begin{prop}
  A Heron tetrahedron can be represented with rational coordinates.
\end{prop}
For the proof I refer to the general result in Section \ref{sec:sn:emb}.\qed

\begin{conj}
  A Heron tetrahedron can be represented with integer coordinates.
\end{conj}
Unfortunately I cannot prove this, but I verified this with the help of a
computer for all the $9019$ Heron tetrahedra ($825$ if multiples are omitted)
with edge length up to $60000$.

It is a special case of Conjecture \ref{conj:2} for $n=3$.

\section{Heron simplices}
\subsection{Basics}
The most interesting thing is how to calculate the volume of a
simplex if only the lengths of the edges are known. It is
possible to generalize the Formula (\ref{eq:v3}) to all
dimensions. One can find a proof in the appendix to my
dissertation \cite{diss}.

\begin{prop}
The volume $V$ of the $n$-dimensional simplex $A_0A_1\cdots A_n$ is given by
\begin{equation}\label{eq:vn}
V^2=\frac{(-1)^{n+1}}{2^n(n!)^2}\det A,
\end{equation}
where $A=(a_{ij})_{i,j=0,\dots,n+1}$ is a $(n+2)\times(n+2)$-matrix with
$a_{ij}=\overline{A_iA_j}^2$ for $i,j\le n$ and $a_{ij}=1-\delta_{ij}$ else.\qed
\end{prop}
\begin{prop}
The $n$-simplex with given lengths of edges exists if and only if all
$(n-1)$-subsimplices (i.e. the ``hyperfaces'' of the simplex) exist and the
right hand side of Equation (\ref{eq:vn}) is positive.\qed
\end{prop}

Now we define: a {\em Heron simplex} is a simplex where all the lengths, areas,
volumes\dots{} are integers.

\subsection{Special types}
First of all: I don't know any example of a Heron simplex with dimension
greater than three. So may be this section discusses the empty set.
So the following seems to be only useful hints for searching Heron simplices.

\subsubsection{Right-angled-vertex simplices}
All face angles in one vertex are right angles. Then the vertices
$A_0,\dots,A_n$ can be represented as $A_0=0$ and $A_i=a_i\cdot {\mathfrak
  e}_i$, $i=1,\dots,n$, where $a_i$ are positive integers and $\{{\mathfrak
  e}_i\}_{i=1,\dots,n}$ is an orthonormal basis.

By the map $x_i=a_i^{-1}$ we get the lengths of a rational $n$-dimensional box 
and vice versa. So the existence problem is equivalent to the existence
problem for a rational box in dimension $n$. But it's not known whether such
boxes exists (see the problem in dimension three).

\subsubsection{Right-angled-face simplices}
All faces are right-angled triangles. Then we have a representation in Cartesian
coordinates $A_0=(0,\dots,0)$, $A_1=(a_1,0,\dots,0)$, $A_2=(a_1,a_2,0,\dots,0)$,
\dots, $A_n=(a_1,a_2,\dots,a_n)$.

The edges have to be integer, so the numbers $\sum\limits_{i=s}^ta_i^2$ 
for $1\le s\le t\le n$ have to be square numbers. All areas, volumes,\dots are 
products of edge lengths divided by an integer depending on the dimension,
so these are at least rationals, but we get integers by scaling.
I don't know whether such sequences with $n>3$ exists or are known.
If one exists, then one of the generators of the Pythagorean triples has to be
greater than $5320$.


\subsection{Embedding properties}\label{sec:sn:emb}
\begin{prop}\label{prop:sn:rat}
  A Heron simplex can be represented with rational coordinates.
\end{prop}
\begin{proof}
I will show a little bit more: it is possible to embedded the simplex
$A_0A_1\cdots A_n$ such that for any $i$ all but the first $i$ coordinates
of $A_i$ vanish.

We do induction by dimension $n$. We only have the step from $n-1$
to $n$ left. Let the subsimplex $A_0A_1\cdots A_{n-1}$ be embedded
as wanted. Let $x_i^{(j)}$ be the $i$th coordinate of $A_j$. Then
the coordinates $x_1,\dots,x_n$ of $A_n$ are determined by
\begin{equation*}
\begin{array}{ccccccccccc}
  x_1^2 & + & x_2^2 & + & x_3^2 & + & \cdots & + & x_n^2 & = &
  \overline{A_0A_n}^2, \\
  (x_1-x_1^{(1)})^2 & + & x_2^2 & + & x_3^2 & + & \cdots & + & x_n^2 & = &
  \overline{A_1A_n}^2, \\
  (x_1-x_1^{(2)})^2 & + & (x_2-x_2^{(2)})^2 & + & x_3^2 & + & \cdots
  & + & x_n^2 & = & \overline{A_2A_n}^2,
\end{array}
\end{equation*}
and so on. Subtracting the first two equations one gets a linear
equation in $x_0$, subtracting the next two using the value of
$x_0$ one gets a linear equation in $x_1$, and so on. So all $x_i$
for $i=0,\dots,n-1$ are rational. The value of $x_n$ is the height
of the simplex, but this is $n$ times the ratio of volume of the
simplex and the volume of the subsimplex $A_0A_1\cdots A_{n-1}$
which are integers.\qed
\end{proof}

\begin{conj}
  A Heron simplex can be represented with integer coordinates.
\end{conj}
By Proposition \ref{prop:sn:rat} this is a special case of Conjecture
\ref{conj:2}.

\section{Embedding Conjecture}
All the presented embedding properties are special cases of the following
conjecture.

\begin{conj}\label{conj:2}
  Let $M\subset \mathbb{Q}^n$ be a set of points such that the distances
  between any two points of $M$ are integer.
  Then one can find a Euclidean motion $T$ such that $TM\subset\mathbb{Z}^n$.

  Short: Integer distances $\Rightarrow$ integer coordinates.
\end{conj}

This conjecture is equivalent to

\begin{conj}\label{conj:1}
  Let $M\subset \mathbb{Z}^n$ be a set of points such that the distances
  between any two points of $M$ are integer and divisible by a number $k$.
  Then one can find a set $N\subset \mathbb{Z}^n$ such that $k\cdot N$ (the
  set $N$ scaled by factor $k$) is congruent to $M$.

  Short: Distances divisible by $k$ $\Rightarrow$ coordinates divisible by $k$.
\end{conj}

\begin{conj}\label{conj:3}
  Let $p$ be prime. Then holds:
  If $M\subset \mathbb{Z}^n$ is a set of points such that the distances
  between any two points of $M$ are integer and divisible by $p$,
  then one can find a Euclidean motion $T$ such that $TM\subset(p\mathbb{Z})^n$.

  Short: Distances divisible by $p$ $\Rightarrow$ coordinates divisible by $p$.
\end{conj}

Obviously for a given number $n$ Conjecture \ref{conj:1} is equivalent to
Conjecture \ref{conj:3} for all primes $p$.

All conjectures are trivial in the case $n=1$.

\begin{prop}\label{prop:em2}
  The Conjecture \ref{conj:3} is true for $n=2$ and all primes $p$,
  i.e. Conjecture \ref{conj:2} is true for $n=2$.
\end{prop}
\begin{proof}
  By translation of $M$ we can assume $(0,0)\in M$. Then holds:
  \begin{equation}
    \label{eq:41}
    (x,y)\in M \qquad\Rightarrow\qquad p^2\mid x^2+y^2.
  \end{equation}
  If $p=2$ or $p\equiv 3(4)$ then immediately follows $p\mid x,y$.

  Let $p\equiv 1(4)$. Then there exists integers $a$ and $b$ satisfying
  $a^2+b^2=p$. In Gaussian numbers Equation (\ref{eq:41}) is equivalent to
  $(a+bi)^2(a-bi)^2\mid (x+yi)(x-yi)$, which means
  \begin{equation}
    \label{eq:42}
    (a+bi)^2\mid x+yi \quad\text{ or }\quad
    p\mid x+yi \quad\text{ or }\quad
    (a-bi)^2\mid x+yi.
  \end{equation}
  We will show that one of $(a+bi)^2$, $p$ and $(a-bi)^2$ is a factor of
  $x+yi$ for all $(x,y)\in M$.
  
  {\em Step 1: One of  $a+bi$ and $a-bi$ is a factor of $x+yi$ for all
  $(x,y)\in M$:}\\
  Assume $a+bi\nmid x_1+y_1i$ and $a-bi\nmid x_2+y_2i$. Since 
  $a-bi\mid x_1+y_1i$ and $a+bi\mid x_2+y_2i$ by Equation (\ref{eq:42}) it
  follows $a\pm bi\nmid (x_1-x_2)+(y_1-y_2)i$ which contradicts
  $p\mid (x_1-x_2)^2+(y_1-y_2)^2$.

  {\em Step 2: Assume $a\pm bi$ is a factor of $x+yi$ for all $(x,y)\in M$. Then
    one of $(a\pm bi)^2$ and $p$ is a factor of $x+yi$ for all
    $(x,y)\in M$:}\\
  Assume $(a\pm bi)^2\nmid x_1+y_1i$ and $p\nmid x_2+y_2i$. From Equation
  (\ref{eq:42}) follows $p\mid x_1+y_1i$ and $(a\pm bi)^2\mid x_2+y_2i$, hence
  $(a\pm bi)^2\nmid (x_1-x_2)+(y_1-y_2)i$ and
  $p\nmid (x_1-x_2)+(y_1-y_2)i$
  which contradicts $a\pm bi\mid (x_1-x_2)+(y_1-y_2)i$ and Equation
  (\ref{eq:42}).

  If $p$ is a factor of all $x+yi$ we are done. If $(a\pm bi)^2$ is a factor
  of all $x+yi$ we set $Tz=\frac{(a\mp bi)^2}{p}\cdot z$ and are done
  ($\left|\frac{(a\mp bi)^2}{p}\right|=1$, so $T$ is a rotation).
  \qed
\end{proof}

\begin{prop}
  The Conjecture \ref{conj:3} is true for $n=3$ and $p\le37$.
\end{prop}
\begin{proof}
  For $p=2$ it is trivial.

  For $3\le p\le37$ it's an easy implication of the following assertion
  which was verified by a computer program:
  \begin{assertion}
    Let $x\in \mathbb{Z}^3$, $p^2\mid \|x\|^2$ and $p\nmid x$. Then there exists
    a matrix $A\in \operatorname{M}(\mathbb{Z},3)$, unique up to
    $\operatorname{O}(\mathbb{Z},3)$, such that $AA^T=p^2\cdot I$ and
    $p^2\mid Ax$. Furthermore for all $y\in \mathbb{Z}^3$ with $p^2\mid \|y\|^2$,
    $p^2\mid\|x-y\|^2$ follows $p^2\mid Ay$.
  \end{assertion}
  \qed
\end{proof}

In dimension $4$ there is a counterexample to Conjecture \ref{conj:3}. We take
the points with coordinates $(0,0,0,0)$, $(1,1,1,1)$, $(2,0,0,0)$ and
$(1,1,1,-1)$ which have pairwise distances of $2$. But it's impossible to
embed a regular tetrahedron of length $1$ into the $\mathbb{Z}^4$. So the
conjectures must be slightly modified for higher dimensions: the points must
not lie in a hyperplane.

\begin{appendix}
  \section{Examples for Heron tetrahedra}
  \subsection{Semi-regular tetrahedra}\label{ap:sr}
  The following table lists all semi-regular tetrahedra with edge
  lengths up to $10000$. The multiples are omitted. The table shows the three
  lengths of the edges and an integer embedding with $A_0=(0,0,0)$.

  {\small
  \begin{tabular}{|ccc|ccc|}
    \hline
    $a$ & $b$ & $c$ & $A_1$ & $A_2$ & $A_3$ \\
    \hline
    203 & 195 & 148 & (168,112,21) & (180,0,-75) & (12,112,-96)\\
    888 & 875 & 533 & (864,192,72) &  (812,-315,-84) &  (164,-123,492)\\
    1804 & 1479 & 1183 & (1452,1056,176) &  (360,1329,540) &  (756,273,-868)\\
    2431 & 2296 & 2175 & (2332,561,396) &  (1792,-1344,-504) &  (1044,-783,1740)\\
    2873 & 2748 & 1825 & (2652,884,663) &  (2652,-576,-432) &  (600,-700,1575)\\
    3111 & 2639 & 2180 & (3111,0,0) &  (1911,1456,1092) &  (1200,448,-1764)\\
    5512 & 5215 & 1887 & (5512,0,0) & (4900,1428,1071) & (612,-1764,-273)\\
    8484 & 6625 & 6409 & (8316,1344,1008) & (5300,-3180,-2385) & (3536,-156,5343)\\
    \hline
  \end{tabular}
  }

  \subsection{Right-angled-face and isosceles-face tetrahedra}\label{ap:raf}
  The following table lists all right-angled-face tetrahedra with edge lengths
  up to $5000$ and the associated isosceles-face tetrahedra. The multiples are
  omitted. Note that the order of $A$, $B$, $C$ and $D$ is not $A_0$, $A_1$,
  $A_2$ and $A_3$.

  {\scriptsize
  \begin{tabular}{|cccccc|ccc|ccc|}
    \hline
    $AB$ & $AC$ & $AD$ & $BC$ & $BD$ & $CD$ &
    $a_1$ & $a_2$ & $a_3$ & $2a$ & $x$ & $2b$ \\
    \hline
    697 & 672 & 680 & 185 & 153 & 104 & 672 & 104 & 153 & 1344 & 697 & 306 \\
    925 & 756 & 765 & 533 & 520 & 117 & 756 & 117 & 520 & 1512 & 925 & 1040 \\
    1073 & 952 & 448 & 495 & 975 & 840 & 448 & 840 & 495 & 896 & 1073 & 990 \\
    1105 & 952 & 1073 & 561 & 264 & 495 & 952 & 495 & 264 & 1904 & 1105 & 528 \\
    1105 & 975 & 1073 & 520 & 264 & 448 & 975 & 448 & 264 & 1950 & 1105 & 528 \\
    2165 & 2040 & 2067 & 725 & 644 & 333 & 2040 & 333 & 644 & 4080 & 2165 & 1288 \\
    2665 & 2175 & 1092 & 1540 & 2431 & 1881 & 1092 & 1881 & 1540 & 2184 & 2665 & 3080 \\
    3277 & 2555 & 1925 & 2052 & 2652 & 1680 & 1925 & 1680 & 2052 & 3850 & 3277 & 4104 \\
    3485 & 2640 & 2652 & 2275 & 2261 & 252 & 2640 & 252 & 2261 & 5280 & 3485 & 4522 \\
    3485 & 2640 & 3179 & 2275 & 1428 & 1771 & 2640 & 1771 & 1428 & 5280 & 3485 & 2856 \\
    3485 & 3360 & 3444 & 925 & 533 & 756 & 3360 & 756 & 533 & 6720 & 3485 & 1066 \\
    3965 & 3723 & 840 & 1364 & 3875 & 3627 & 840 & 3627 & 1364 & 1680 & 3965 & 2728 \\
    4181 & 3740 & 4100 & 1869 & 819 & 1680 & 3740 & 1680 & 819 & 7480 & 4181 & 1638 \\
    4225 & 4180 & 4199 & 615 & 468 & 399 & 4180 & 399 & 468 & 8360 & 4225 & 936 \\
    4453 & 3485 & 2275 & 2772 & 3828 & 2640 & 2275 & 2640 & 2772 & 4550 & 4453 & 5544 \\
    \hline
  \end{tabular}
  }
  
  \subsection{Generic tetrahedra}\label{ap:gen}
  The following table lists all generic (i.e. not contained in the previous sections)
  tetrahedra with edge lengths up to $1000$. The multiples are omitted. The
  table shows the lengths of the edges $A_0A_1$, $A_0A_2$, $A_0A_3$,
  $A_1A_2$, $A_1A_3$, $A_2A_3$, and an integer embedding with $A_0=(0,0,0)$.

  {\footnotesize
    \begin{tabular}{|cccccc|ccc|}
      \hline
      \multicolumn{6}{|c|}{edges} &
      $A_1$ & $A_2$ & $A_3$\\
      \hline
      117 & 84 & 80 & 51 & 53 & 52 &
      (108,36,27) & (84,0,0) & (64,48,0) \\
      160 & 153 & 120 & 25 & 56 & 39 &
      (128,96,0) & (108,108,9) & (72,96,0) \\
      225 & 200 & 87 & 65 & 156 & 119 &
      (180,108,81) & (120,128,96) & (36,72,33) \\
      318 & 221 & 221 & 203 & 175 & 42 &
      (288,126,48) & (176,-21,132) & (176,21,132) \\
      319 & 318 & 221 & 175 & 210 & 175 &
      (231,176,132) & (126,288,48) & (21,176,132) \\
      319 & 318 & 221 & 175 & 252 & 203 &
      (231,176,132) & (126,288,48) & (-21,176,132) \\
      429 & 300 & 176 & 261 & 275 & 140 &
      (396,132,99) & (288,-84,0) & (176,0,0) \\
      468 & 340 & 297 & 232 & 225 & 65 &
      (432,144,108) & (336,-48,20) & (297,0,0) \\
      595 & 429 & 325 & 208 & 276 & 116 &
      (588,84,35) & (396,132,99) & (312,84,35) \\
      595 & 507 & 325 & 116 & 276 & 208 &
      (588,91,0) & (504,27,48) & (312,91,0) \\
      595 & 555 & 429 & 100 & 208 & 204 &
      (588,84,35) & (540,120,-45) & (396,132,99) \\
      612 & 455 & 480 & 319 & 156 & 185 &
      (432,432,36) & (399,168,-140) & (384,288,0) \\
      671 & 663 & 225 & 580 & 544 & 444 &
      (528,396,121) & (468,204,-423) & (144,108,-135) \\
      680 & 615 & 672 & 185 & 104 & 153 &
      \multicolumn{3}{|c|}
      {edge perpendicular to face\textsuperscript{*}, pair 1} \\
      680 & 680 & 615 & 208 & 185 & 185 &
      \multicolumn{3}{|c|}
      {isosceles faces with common basis\textsuperscript{*}, pair 1} \\
      697 & 697 & 672 & 306 & 185 & 185 &
      \multicolumn{3}{|c|}
      {$2\times$ right-angled-face tetrahedron\textsuperscript{*}} \\
      741 & 477 & 475 & 420 & 304 & 340 &
      (684,228,171) & (432,-108,171) & (380,228,171) \\
      780 & 765 & 715 & 219 & 221 & 148 &
      (624,468,0) & (648,360,189) & (660,264,77) \\
      884 & 880 & 715 & 84 & 205 & 187 &
      \multicolumn{3}{|c|}
      {edge perpendicular to face\textsuperscript{*}, pair 2} \\
      884 & 884 & 715 & 168 & 205 & 205 &
      \multicolumn{3}{|c|}
      {isosceles faces with common basis\textsuperscript{*}, pair 2} \\
      935 & 928 & 900 & 345 & 145 & 260 &
      (660,561,352) & (384,768,352) & (540,576,432) \\
      990 & 795 & 598 & 663 & 436 & 427 &
      (792,594,0) & (720,27,336) & (552,230,0) \\
      990 & 901 & 793 & 793 & 901 & 308 &
      (792,594,0) & (396,451,672) & (396,143,672) \\
      \hline
    \end{tabular}
    }\\
  \textsuperscript{*}: The embedding can be obtained by the construction of
  these tetrahedra, see Figure \ref{fig:ra2if} and the remarks at the end of
  Section \ref{sec:othertypes}.
  
\end{appendix}

\begin{tabular}{ll}
  Authors address:&
  Jan Fricke \\
  & Institut f\"ur Mathematik und Informatik\\
  & Ernst-Moritz-Arndt Universit\"at Greifswald\\
  & Jahnstr. 15a\\
  & D-17487 Greifswald, Germany\\[1em]
  E-Mail:&
  \href{mailto:fricke@uni-greifswald.de}{\tt fricke{\rm @}uni-greifswald.de}
\end{tabular}

\newpage

\tableofcontents

\end{document}